\documentstyle[12pt,epsf]{article}

\begin{document}

\title {Is it possible to estimate the oscillating sum
{\large
$Z=e^{-N}\sum\limits_{n=0}^\infty\frac{N^n}{n!}
\exp
\biggl\{i\biggl(\frac{An}{\sqrt N}+\frac{Bn^2}{2\sqrt{N^3}}\biggr)\biggr\}$
}
 for {\large $N=10^{23}\;$} ?
\author{Alexander~M.~Chebotarev
\thanks{This work was partially supported by
INTAS-99-0698
}\\
{\small 119899 Moscow, Moscow State University, Quantum Statistics Dept.}}}
\date {}
\maketitle

\begin{abstract}
We prove that
it is possible to approximate $Z$
by an  integral $G$
with respect to the  standard Gaussian measure
such that the approximation error
$|G-Z|=O(N^{-1})$
is small,
and the value of the integral
$$
Z=\frac{1}{\sqrt{D}}
e^{i\sqrt{N}(A+B/2)-(A+B)^2/2D}
\biggl(1-i\frac{(A+B)^3}{6\sqrt{N}}\biggr)+
O(N^{-1})
$$ with $D=(1-iB/\sqrt{N})$, is finite
and oscillates rapidly in $A$ and $B$.
Straightforward summation of $Z$  is impossible because
over the huge
interval $n\in[N-\sqrt{2N\log N},N+\sqrt{2N\log N}]$, the summands
contribute significantly.
We discuss approximation for some oscillating  series related to $Z$.
\end{abstract}

\section{Introduction}
When a mathematician says "consider the asymptotics of
a state $\psi_\hbar(x)$ as
the Plank constant $\hbar$ goes to zero", this makes no
rigorous sense
in physics, because $\hbar$ is a constant. The difference between
the mathematical infinity and huge physical
 quantities such as the Avogadro number $N_A$ or $\hbar^{-1}$ is that
$\log N_a\approx 23$, and $\log \hbar^{-1}\approx 37$  are
rather moderate
numbers, if they are used as multipliers. Therefore, for mathematicians,
$1/{\log N_A}=o(1)$ and $1/{\log \hbar^{-1}}=o(1)$
(if they
keep $N_A\to\infty$ and $1/\hbar\to\infty$ in mind),
and in physics,
$1/{\log N_A}=O(1)$, $1/{\log \hbar^{-1}}=O(1)$,
${\log N_A}/{N_A}=O(N_A^{-1})$ and ${\hbar}/{\log \hbar^{-1}}=O(\hbar)$,
e.~t.~c.   Fortunately,
a rather general property of the nature is
that the logarithms of overwhelming physical quantities
are {\sl  moderate numbers}. On the other hand,
the logarithms of large parameters
occur quite rare in asymptotic expansions. But we will see that
 crucial estimates of the asymptotic expansion of $Z$
relay on the fact that $\log N$ is finite.

The sum $Z$ represents some quantities related to explicit solution of
a problem on interoferometrical detection of
gravitational waves
(see \cite{BJK}--\cite{Ch}),
where $N$ is the average number of photons generated
by a laser, and $A=A(t), \;B=B(t)$ are small time-dependent functions
describing interferation of photons splitted into two arms of the
measuring device.
The main computational difficulty is that $2\sqrt{2N\log N}\approx
10^{12}$ summands
 centered at $n=N=10^{23}$ give a relevant
contribution to $Z$. Thus the straightforward computation
of $Z$ is impossible.

In this paper we prove that the series $Z(A,B,N)$
can be approximated by the Gaussian integral, because
in terms of the new variable $x_n=\frac{n-N}{\sqrt{N}}$,
$dx_n=\frac{1}{\sqrt N}$
($n=0,1,\dots$), the Stirling expansion  of Gamma--function implies
the Gaussian approximation of the Poisson probability
distribution $P_N(n)$:
$$
P_N(n)\stackrel{def}{=}e^{-N}\frac{N^n}{n!}=\frac{e^{
-\frac{x_n^2}{2}}}{\sqrt{2\pi N}}
\frac{e^{\frac{x_n^3}{6\sqrt{N}}-\frac{x_n^4}{12N}}}
                       {
\biggl(1+\frac{x_n}{2\sqrt{N}}+\frac{1}{12N}\biggr)
}
\bigl(1+O\bigl(N)^{-3/2}\bigr)
\eqno(1.1)
$$
if $x_n/\sqrt{N}$ is small \cite{Fe} (see Appendix 3, (A.6), (A.7)).

\begin{figure}[ht]
\epsfxsize=8cm
\centerline{\epsffile{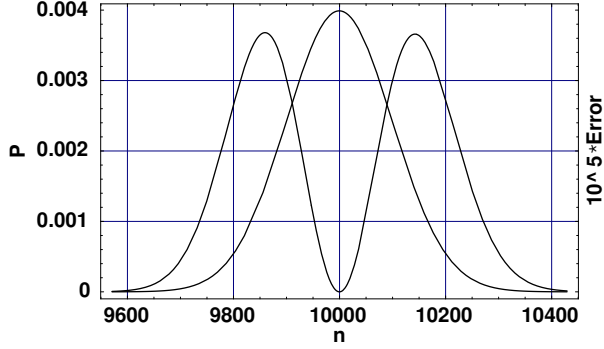}}
\caption{The picture shows the Poisson distribution $P_N(n)$
 for $N=10^4$,
which practically coincides with the Gaussian distribution $N(10^4,10^4)$,
and the approximation error (double-well curve scaled by
$10^5$), representing the
difference between the Poisson and normal
distributions with the Stirling corrections.
}
\label{fig1}
\end{figure}

On the other hand, by the trapezoidal approximation of the integral,
it is possible to estimate the sum
$$
Z=\sum_n P_N(n)e^{iS_N(n)},
\quad S_N(n)=\frac{An}{\sqrt{N}}+\frac{Bn^2}{2\sqrt{N^3}}
$$
over the interval
$B_N\stackrel{def}{=}[N-\sqrt{2N\log N},N+\sqrt{2N\log N}]\subset
{\bf N}.$
In this case, $\max_{n\in B_N} |x_n|=\sqrt{2\log N}$.
More precise, to estimate the sum $Z$, we must
justify several approximations:
\begin{itemize}
\item[-]  estimate the error of
the restriction to $B_N$ of the sum $Z$;
\item[-]  estimate  on $B_N$ the  error of approximation
of the Poisson distribution $P_N(b)$ by the Gaussian one as in Eq. (1.1);
\item[-]  estimate the approximation error of
the replacing of the sum over $B_N$ by the Gaussian integral over the set
$$
D_N\stackrel{def}{=}[-\sqrt{2\log N},\,\sqrt{2\log N}]\subset {\bf R}
$$
\item[-]  expand (1.1) in $1/\sqrt{N}$ on $D_N$
and estimate the approximation error;
\item[-]  extend the Gaussian integral to $\bf R$ and
evaluate it analytically.
\end{itemize}

First, in Section 2, we cut out the left and right tails of the probability distribution
$P_N(n)$ outside of $B_N$, and then in Section 3,
we approximate the sum by an integral
over the interval $D_N\subset\bf R$
corresponding to $B_N$:
$$
Z=
\frac{1}{2\sqrt{N}}
\sum_{x_n\in B_N}\biggl(f_N(x_n)e^{iS_N(n)}
+f_N(x_{n+1})e^{iS_N(n+1)}\biggr)+\eta_N=
$$
$$
=\int_{x\in D_N}\,
dx\,f_N(x)e^{i\sigma_N\bigl(x+\frac{1}{2\sqrt{N}}\bigr)}+
\eta_N+\epsilon_N,
\eqno(1.2)
$$
where
$\epsilon_N$ and $\eta_N$  are the approximation errors,
$\sigma_N(x)=S_N(x\sqrt{N}+N)$
and $f_N(x)=\sqrt{N}\,P_N(x\sqrt{N}+N)$.
More precise, for the total approximation error
$\eta_N+\epsilon_N=Z-\int\,dx\,f_N(x)e^{i\sigma_N(x)}$, we have
$$
|\eta_N|+|\epsilon_N|\le
\sum_{n\in {\bf N}\setminus B_N }P_N(n)+
$$
$$
+\frac{1}{6}\sum_{n}
\biggl\{ \max_{y\in I_n}|P_N^{(2)}(y)|
+2 \max_{y\in I_n}|P_N^{(1)}(y)S_N^{(1)}(y)| +
\max_{y\in I_n}|P_N(y)S_N^{(2)}(y)|\biggr\},
\eqno(1.3)
$$
where $I_n\stackrel{def}{=}[n,n+1]$  (see Appendix (A.1) and (A.5)).

In what follows we prove that
 $|\epsilon_N|+|\eta_N| =O(N^{-1})$, and Eq.~(1.1) describe
 convenient analytical approximation of $f_N(x)$ on $D_N$ by the perturbed
 Gaussian distribution (see (A.6) and (A.7)).

 This implies that
the sum $Z$ can be approximated by the Gaussian integral
$$
Z=\frac{1}{\sqrt{2\pi}}\int_{x\in D_N}\,dx\,
e^{-\frac{x^2}{2}}
\frac{e^{\frac{x^3}{6\sqrt{N}}-\frac{x^4}{12N}}}
{\biggl(1+\frac{x}{2\sqrt{N}}+\frac{1}{12N}\biggr)}
e^{iS_N(\sqrt{N}x+N)}+
O(N^{-1}).
$$
Eventually, the last integral can be extended to ${\bf R}$ and
evaluated explicitly up to terms of order
$O(N^{-1})$:
$$
Z=
\frac{e^{i(A+B/2)\sqrt{N}}}{\sqrt{1-iB/\sqrt{N}}}
e^{-\frac{(A+B)^2}{2(1-iB/\sqrt{N})}}
\biggl(1-i
\frac{(A+B)^3}{6\sqrt{N}}\biggr)+
O(N^{-1}).
\eqno(1.4)
$$
The last estimate of the error is proved by a generalization of
the Komatsu inequality (3.3) proved in Appendix 4.
This asymptotic expansion  (1.4) is the main result of the note.
It can be applied in many different ways.

In Section 4, we present similar estimates for several related oscillating
series. First, it justifies the asymptotic expansion of the sum
$$
Z_F=e^{-N}\sum\limits_{n=0}^\infty\frac{N^n}{n!}
F\biggl(\frac{n}{\sqrt N}\biggr)
e^{i\frac{qn^2}{2N^{3/2}}},
$$
where $F(x)$ is the Fourier transform of any finitely supported
function $\widetilde F(p)$ such that
$({\rm diam\,supp\,}\widetilde F+|B|)^2<\log N$:
$$
Z_F\approx\frac{1}{\sqrt{2\pi}}\int_R dp\,\widetilde F(p)
\frac{e^{i(p+q/2)\sqrt{N}}}{\sqrt{1-iq/\sqrt{N}}}
e^{-\frac{(p+q)^2}{2(1-iq/\sqrt{N})}}
\biggl(1-i
\frac{(p+q)^3}{6\sqrt{N}}\biggr).
$$

It can also be applied for the correlation and Fourier analysis of
$Z(t)$ in the case, when $A=A(t)$ and $B=B(t)$.

Finally, we note that the series $Z(x,t, \hbar^{-2})
\stackrel{def}{=}\psi_\hbar(x,t)$
represents an explicit solution
of the Cauchy problem for the Schr\"odinger
equation with a special  periodic initial condition:
$$
i\hbar\partial_t\psi=\frac{\hbar^2}{2}\partial^2_x\psi,\quad
\psi|_{t=0}=
\exp\bigl\{-\hbar^{-2}\bigl(1-e^{ix/\hbar}\bigr)\bigr\}.
$$

\section{Estimate of the weight of tails}

The total weight $T_N$ of the left and right tails of the
probability distribution
$P_N(n)$ on the complement of the interval
$B_N$
can be estimated by use of the expansion (1.1) which is clearly valid
on this set. By definition, we have
$$
T_N\stackrel{def}{=}1-\sum_{n\in B_N}P_N(n)=
1-\sum_{x_n\in D_N}
\frac{1}{\sqrt{N}}f_N(x_n).
\eqno(2.1)
$$
The trapezoidal approximation (see Appendix 1, (A.1) for the estimates of
the approximation errors) implies the following
integral representation of the sum:
$$
\sum_{x_n\in D_N}
\frac{1}{\sqrt{N}}f_N(x_n)=\int_{D_N}
dx\, f_N(x) +\gamma_N,
$$
where $|\gamma_N|\le
\frac{1}{6}\sum_{n\in B_N}\max_{x\in I_n}|P^{(2)}_N(x)|$.
Note that the function $P_N(x)$ has just one extremum on the set
$B_N$. Hence its first and second derivatives may have no more then two
and three critical points respectively. Hence the function
$|P^{(2)}_N(x)|$ is monotone on all intervals $I_n$ with the exception
of the three (see Fig.  2).

\begin{figure}[ht]
\epsfxsize=8cm
\centerline{\epsffile{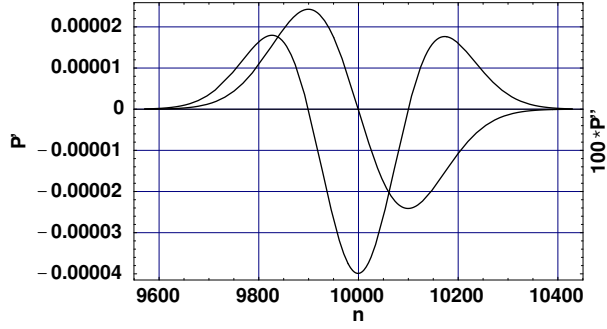}}
\caption{The derivatives of the Poisson distribution
become more flatten as $N$ increases. Each derivation add a new
critical point. The picture shows $P^{(1)}_N(n)$  with  two
critical points, and
$100*P^{(2)}_N(n)$ with three critical points, for $N=10^4$.}
\label{fig2}
\end{figure}

Therefore, we can apply the estimate (A.5):
$$
|\gamma_N|\le \sum_{n\in B_N}\max_{x\in I_n}|P^{(2)}_N(x)|\le 3
\sum_{n\in B_N}|P^{(2)}_N(n)|= O(N^{-1}).
$$
Now the expansion (1.1) ensures the equality:
$$
\int_{D_N}
dx\, f_N(x)=\frac{1}{\sqrt{2\pi}}
\int_{D_N}
dx\,e^{-\frac{x^2}{2}}
\frac{e^{\frac{x^3}{6\sqrt{N}}-\frac{x^4}{12N}}}
                       {
\biggl(1+\frac{x}{2\sqrt{N}}+\frac{1}{12N}\biggr)}=
$$
$$
=\frac{1}{\sqrt{2\pi}}\int_{R}
dx\,e^{-\frac{x^2}{2}}\biggl( 1-\frac{3x-x^3}{6\sqrt{N}}
 \biggr)+O(N^{-1})=1+O(N^{-1})
$$
because the functions $x$ and $x^3$ are odd.
The estimate of the approximation error follows from
the Komatsu inequality \cite{IM} (see Appendix 4):
$$
4e^{-\frac{a^2}{2}}\biggl(\sqrt{a^2+4}+a  \biggr)^{-1}\le
\int_{|x|\ge a} dx\, e^{-\frac{x^2}{2}}\le
4e^{-\frac{a^2}{2}}\biggl(\sqrt{a^2+2}+a  \biggr)^{-1}.
$$
For $a=\sqrt{2\log N}$, we obtain
$$
\frac{1}{\sqrt{2\pi}}\int_{R\setminus D_N}
dx\,e^{-\frac{x^2}{2}}\biggl( 1-\frac{3x-x^3}{6\sqrt{N}}\biggr)=
\frac{1}{\sqrt{2\pi}}\int_{R\setminus D_N}dx\,e^{-\frac{x^2}{2}}\le
$$
$$
\le\frac{4}{\sqrt{2\pi}}
e^{-\frac{a^2}{2}}\biggl(\sqrt{a^2+2}+a  \biggr)^{-1}
\biggr|_{a=\sqrt{2\log N}}=O(N^{-1}).
$$
These estimates explain the specific choice of the sets $D_N$
and $B_N$.
Therefore, the total weight $T_N$ of the tails on the complement of
$B_N$ equals $O(N^{-1})$:
$$
T_N=1-\sum_{n\in B_N}P_N(n)=
1-\sum_{x_n\in D_N}
\frac{1}{\sqrt{N}}f_N(x_n)=
$$
$$
=1-
\frac{1}{\sqrt{2\pi}}
\int_{D_N}
dx\,e^{-\frac{x^2}{2}}
\frac{e^{\frac{x^3}{6\sqrt{N}}-\frac{x^4}{12N}}}
                       {
\biggl(1+\frac{x}{2\sqrt{N}}+\frac{1}{12N}\biggr)}
-\gamma_N\le O(N^{-1})+|\gamma_N|=O(N^{-1}).
\eqno(2.2)
$$

\section{Estimates of the oscillating sum}

The main property of the set $B_N$ is that the Stirling
formula holds true on this set. On the other hand,
the probability distribution $P_N(n)$ is supported on this
set, up to events with total probability of order $O(N^{-1})$.
Hence the expectation of the oscillating exponents can be restricted
to $B_N$ with the same approximation error:
$$
Z=\sum_{n\in B_N}P_N(n)e^{iS_N(n)} +O(N^{-1}).
\eqno(3.1)
$$
The trapezoidal approximation method is applicable to the sum (3.1)
with the approximation error (1.3) (see (A.1) for details):
$$
Z=\int_{x\in D_N} dx\,f_N(x)e^{i\sigma_N(x)} +O(N^{-1})+\gamma_N,
$$
where $$|\gamma_N|\le
\frac{1}{6}\sum_{n\in B_N}
\biggl\{ \max_{y\in I_n}|P_N^{(2)}(y)|
+2 \max_{y\in I_n}|P_N^{(1)}(y)S_N^{(1)}(y)| +
\max_{y\in I_n}|P_N(y)S_N^{(2)}(y)|\biggr\}.$$
By definition of the set $D_N$, we have the following
uniform estimates of the
derivatives
$$
S_N^{(1)} =\frac{A}{\sqrt{N}}+\frac{Bx}{\sqrt{N^3}}=O(N^{-1/2}),\quad
S_N^{(2)} =\frac{B}{\sqrt{N^3}}=O(N^{-3/2}),
\quad x\in D_N
$$
and as is proved in Appendix (see (A.3) and (A.4)),
$$
P_N^{(1)}(x)=P_N(x)\biggl( \log\frac{N}{x}+O(x^{-1}) \biggr)=
P_N(x)\,O(N^{-1/2}),
$$
$$
P_N^{(2)}(x)=P_N(x)\biggl( \biggl(\log\frac{N}{x}\biggr)^2+
x^{-1}\log\biggl(\frac{N}{x}\biggr)
+O(x^{-1}\log x) \biggr)=
P_N(x)\cdot O(N^{-1}),
$$
because $N/x=1+O(N^{-1}\log N)^{1/2}$ uniformly on $D_N$, and hence
$\log(N/x)=O(N^{-1/2})$ and $(\log(N/x))^2=O(N^{-1})$ uniformly
in $x\in D_N$.

Since $\sum_n P_N(n)=1$, it follows that
$|\gamma_N|=O(N^{-1})$. Thus we proved that
$$
Z=
\frac{1}{\sqrt{2\pi}}\int_{x\in D_N}\,dx\,
e^{-\frac{x^2}{2}}
\frac{e^{\frac{x^3}{6\sqrt{N}}-\frac{x^4}{12N}}}
{\biggl(1+\frac{x}{2\sqrt{N}}+\frac{1}{12N}\biggr)}
e^{iS_N(\sqrt{N}x+N)}+
O(N^{-1}),
$$
where $S_N(\sqrt{N}x+N)=A(x+\sqrt{N})+B(x+\sqrt{N})^2/(2\sqrt{N})$.
Now we expand the exponent into the series in $x/\sqrt{N}$, and
extend the  of integration domain to the whole real line:
$$
Z=\frac{1}{\sqrt{2\pi}}\int_{R}\,dx\,
e^{-\frac{x^2}{2}}
\biggl(1-\frac{3x-x^3}{6\sqrt{N}}\biggr)
e^{iA(x+\sqrt{N})+iB(x+\sqrt{N})^2/(2\sqrt{N})}+
O(N^{-1}).
\eqno(3.2)
$$

The error $\sigma_N$ related to the extension of the Gaussian
integral from $D_N$ to ${\bf R}$ can be estimated by an inequality
of Komatsu type:
$$
\int_{|x|\ge a} dx\,x^n e^{-\frac{x^2}{2}}\le
4(a+1)^{n-1}e^{-\frac{a^2}{2}}\biggl(\sqrt{a^2+2}+a  \biggr)^{-1}
\eqno(3.3)
$$
(see Appendix 4), which implies the required absolute estimate
$$
\frac{1}{\sqrt{2\pi}}\int_{|x|\in D_N}\,dx\,x^ne^{-ax^2/2}
\biggr|_{a=\sqrt{2\log N}}\le (\log N)^{n/2}\cdot O(N^{-1})=
O(N^{-1}),
$$
$n=1,2$.

\begin{figure}[ht]
\epsfxsize=8cm
\centerline{\epsffile{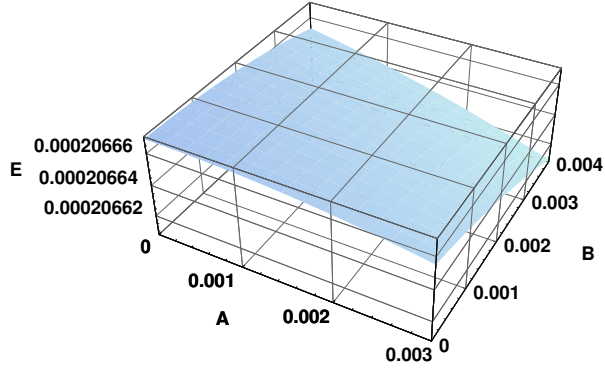}}
\caption{The picture shows the absolute approximation error $E$
of the sum $Z$ by integral (3.2), as a decreasing function of $A$ and $B$
for $N=1000$.}
\label{fig3}
\end{figure}

By using the standard technic, we find elementary Gaussian
integral (3.2):
$$
\frac{1}{\sqrt{2\pi}}\int_{\bf R}\,dx\,xe^{-ax^2/2+ibx}=
\frac{ib}{a^{3/2}}e^{-b^2/2a},
$$
$$
\frac{1}{\sqrt{2\pi}}\int_{\bf R}\,dx\,x^3e^{-ax^2/2+ibx}=
\frac{ib(3a-b^2)}{a^{7/2}}e^{-b^2/2a}.
$$
Hence, by substituting $a=1-iB/\sqrt{N}$, $b=A+B$, we obtain

\medskip\noindent
{\bf Theorem 3.1}
$$
Z=\frac{e^{i\sqrt{N}(A+B/2)}}{\sqrt{1-iB/\sqrt{N}}}
e^{-\frac{(A+B)^2}{2(1-iB/\sqrt{N})}}
\biggl(1-i
\frac{(A+B)^3}{6\sqrt{N}}\biggr)+
O(N^{-1}).
\eqno(3.4)
$$

\section{Approximation of related series}
In this section we consider the integral approximation of the related series
$$
Z_s=e^{-N}\sum\limits_{n=1}^\infty\frac{N^n}{n!}n^s
\exp
\biggl\{i\biggl(\frac{An}{N^{1/2}}+\frac{Bn^2}{2N^{3/2}}\biggr)\biggr\}
\eqno(4.1)
$$
for $s=1,2$, which are used to calculate the photo current and its dispersion
(see \cite{Ch}, Chap. 2). Note that $Z_1$ can be rewritten as
$$
Z_1=Ne^{-N}\sum\limits_{n=0}^\infty\frac{N^n}{n!}
\exp
\biggl\{i\biggl(\frac{A(n+1)}{N^{1/2}}+
\frac{B(n+1)^2}{2N^{3/2}}\biggr)\biggr\}=
$$
$$
=N e^{i\frac{A}{\sqrt{N}}+i\frac{B}{2\sqrt{N^3}}}\cdot Z(A+B/N,B,N)=
N e^{i\frac{A}{\sqrt{N}}}\cdot Z(A,B,N)+O(1).
$$
This relation readily implies
$$
Z_1=N e^{i\frac{A}{\sqrt{N}}}
\frac{e^{i(A+B/2)\sqrt{N}}}{\sqrt{1-iB/\sqrt{N}}}
e^{-\frac{(A+B)^2}{2(1-iB/\sqrt{N})}}
\biggl(1-i
\frac{(A+B)^3}{6\sqrt{N}}\biggr)+
O(1).
$$
Similarly,
$$
Z_2=Z_1+
N^2e^{-N}\sum\limits_{n=0}^\infty\frac{N^n}{n!}
\exp
\biggl\{i\biggl(\frac{A(n+2)}{N^{1/2}}
+\frac{B(n+2)^2}{2N^{3/2}}\biggr)\biggr\}=
$$
$$
=Z_1+N^2 e^{i\frac{2}{\sqrt{N}}(A+B/N)}\cdot Z(A+2B/N, B,N)
+O(N)=
$$
$$
=N^2 e^{i\frac{2A}{\sqrt{N}}}Z(A,B,N)+O(N).
\eqno(4.2)
$$
For the general case, we can prove the following assertion.

\medskip\noindent
{\bf Theorem 4.1} {\sl  The sum} (4.1) {\it has the following
asymptotic expansion}
$$
Z_s=N^s e^{i\frac{sA}{\sqrt{N}}}Z(A,B,N)+O(N^{s-1}),
$$
{\it where $Z(A,B,N)$ is given by Eq.}~(3.4).
\medskip

Another series which can be estimated with the help of previously
developped technics is
$$
\widetilde Z=e^{-N}\sum\limits_{n=0}^\infty\frac{N^n}{n!}
\exp
\biggl\{i\biggl(A\sqrt n+\frac{Bn}{2\sqrt N}\biggr)\biggr\}
\eqno(4.3)
$$
Indeed, for the new variable $x=(n-N)/\sqrt{N}$, we have
$n/\sqrt{N}=x+\sqrt{N}$, and
$$
\sqrt{n}=(x+\sqrt{N})(1+x/\sqrt{N})^{-1/2}=
\sqrt{N}+\frac{x}{2}-\frac{3x^2}{8\sqrt{N}}+O(N^{-1}).
$$
Hence,
$$
A\sqrt{n}+\frac{Bn}{2\sqrt{N}}=
A\biggl(\sqrt{N}+\frac{x}{2}-\frac{3x^2}{8\sqrt{N}}\biggr)+
B(\sqrt{N}+x)/2,
\eqno(4.4)
$$
and by replacing the exponent in Eq.~(3.2) by (4.4),
 we obtain  the following assertion

\medskip\noindent
{\bf Theorem 4.2} {\it The sum} (4.3) {\it has the following asymptopic
exapnsion}
$$
\widetilde Z=
\frac{e^{i\sqrt{N}(A+B/2)}}{\sqrt{1+3iA/4\sqrt{N}}}
e^{-\frac{(A+B)^2}{8(1+3iA/4\sqrt{N})}}
\biggl(1-i
\frac{(A+B)^3}{48\sqrt{N}}\biggr)+
O(N^{-1}).
\eqno(4.5)
$$

To conclude this section, we consider an estimate of
the following multiple sum:
$$
{\bf Z}=e^{-2N}\sum\limits_{m,n=0}^\infty\frac{N^{m+n}}{m!\,n!}
\exp
\biggl\{i\biggl(\frac{a_1m+a_2n}{\sqrt{N}}
+\frac{b_1m+b_1n+2b_3mn}{2\sqrt{N^3}}\biggr)\biggr\}
\eqno(4.6)
$$
which was derived in \cite{BM} to describe the interaction of the lazer beam
in some standard state with the quantum oscillator.

In variables $x=(m-N)/\sqrt N$, $y=(n-N)/\sqrt N$, this series can be
approximated by the Gaussian integral
$$
{\bf Z}=\frac{e^{i\sqrt N \sigma(a,b)}}{2\pi}
\int_{{\bf R}^2}\,dx\,dy\,e^{S_N(x,y)}
\biggl(1-\frac{3(x+y)-x^3-y^3}{6\sqrt N}\biggr)+O(N^{-1})
$$
where $\sigma(a,b)=a_1+a_2+b_3+(b_1+b_2)/2$ and
$$
S_N(x,y)=-(x^2+y^2)/2+
$$
$$
+i(a_1+b_1+b_3)x+i(a_2+b_2+b_3)y+i(b_1x^2+b_2y^2+2b_3xy)/2\sqrt N.
$$
This integral can easily be evaluated. As a result we obtain the following
estimate.

\medskip\noindent
{\bf Theorem 4.3} {\it The sum} (4.6) {\it has the following asymptopic
exapnsion}
$$
{\bf Z}=\frac{e^{i\sqrt N \sigma(a,b)- \Theta_N(a,b)}}
{\sqrt{1-i(b_1+b_2)/\sqrt N}}\biggl(1-i
\frac{(a_1+b_1+b_3)^3+(a_2+b_2+b_3)^3}{6\sqrt{N}}\biggr)+
O(N^{-1}),
$$
{\it where}
$$
\Theta_N(a,b)=
\frac{(a_1+b_1+b_3)^2+(a_2+b_2+b_3)^2}{2{\sqrt{1-i(b_1+b_2)/\sqrt N}}}
+i\frac{b_3(a_1+b_1+b_3)(a_2+b_2+b_3)}{\sqrt N}.
$$

\section{Appendix: Simple a priori estimates}
{\small

\subsection{Trapezoidal approximation}
Consider the estimate of the trapezoidal approximation of  the integral of
a smooth function $f(x)$.
The Taylor expansion implies
$$
\frac{1}{2}\bigl(f(x_n)+f(x_{n+1})\bigr)=
f(x_{n+\frac{1}{2}})+\kappa_n+\kappa_{n+\frac{1}{2}},
$$
where  the approximation error does not exceed the value
$$
|\kappa_n|\le \frac{1}{2}(x_n-x_{n+\frac{1}{2}})^2
\max_{y\in [x_n,x_{n+\frac{1}{2}}]}|f^{(2)}(y)|,
$$
and for $x\in[x_n,x_{n+1}]$, we have
$$
f(x)=f(x_{n+\frac{1}{2}})+(x-x_{n+\frac{1}{2}})f^{(1)}(x_{n+\frac{1}{2}})
+\gamma_n(x),
$$
where $|\gamma_n(x)|\le \frac{1}{2}(x-x_{n+\frac{1}{2}})^2
\max_{y\in [x_n,x_{n+1}]}|f^{(2)}(y)|$.
Hence the error of the trapezoidal approximation can be estimated
as follows
$$
\int_{x_n}^{x_{n+1}}f(x)\,dx=
\frac{x_{n+1}-x_n}{2}\biggl(f_n+f_{n+1}-
\kappa_n-\kappa_{n+\frac{1}{2}}\biggr)-\int_{x_n}^{x_{n+1}}\gamma_n(x)\,dx=
$$
$$
=\frac{x_{n+1}-x_n}{2}\bigl(f_n+f_{n+1}\bigr)+\eta_n,
\eqno(A.1)
$$
where $f_n=f(x_n)$ and $|\eta_n|\le \frac{1}{6}(x_{n+1}-x_n)^3
\max_{y\in [x_n,x_{n+1}]}|f^{(2)}(y)|$.

\medskip
\subsection{Derivatives of the Poisson distribution}
Consider the estimates of the derivatives of the distribution
$P_N(x)$. Note that from the Stirling expansion for the Gamma--function
$$
\Gamma(x+1)=\biggl(\frac{x}{e}\biggr)^x\sqrt{2\pi x}\biggl( 1+\frac{1}{12x}
+\frac{1}{288x^2}+O(x^{-3})\biggr)
\eqno(A.2)
$$
we have the following estimate for its derivatives:
$$
\Gamma^{(1)}(x+1)=
\biggl(\log x +\frac{1}{2x}+O(x^{-2})  \biggr)\Gamma(x+1)
$$
and similarly,
$$
\Gamma^{(2)}(x+1)=\biggl\{
\biggl(\log x+\frac{1}{2x}+O(x^{-2})  \biggr)^2+
\biggl( \frac{1}{x}+O(x^{-2}) \biggr)\biggr\}\Gamma(x+1)=
$$
$$
\biggl\{\log^2x+\frac{\log x}{x}
+\frac{1}{x}+O(x^{-2}\,\log x)\biggr\}\Gamma(x+1).
$$
Therefore,
$$
\biggl(\frac{\Gamma^{(1)}(x+1)}{\Gamma(x+1)}\biggr)^2-
\frac{\Gamma^{(2)}(x+1)}{\Gamma(x+1)}=O(x^{-1}\log x),
\eqno(A.3)
$$
and moreover,
$$
P^{(1)}_N(x)=e^{-N}\frac{d}{dx}\frac{N^x}{\Gamma(X+1)}=
P_N(x)\biggl( \log N- \frac{\Gamma^{(1)}(x+1)}{\Gamma(x+1)}\biggr)=
$$
$$
= P_N(x)\biggl( \log \frac{N}{x}-\frac{1}{2x}+O(x^{-2})\biggr).
\eqno(A.4)
$$
Similarly, for the second derivative we have
$$
P^{(2)}_N(x)=e^{-N}\frac{d^2}{dx^2}\frac{N^x}{\Gamma(X+1)}=
$$
$$
=P_N(x)\biggl\{\biggl(\log \frac{N}{x}-\frac{1}{2x}+O(x^{-2})\biggr)^2-
\frac{\Gamma^{(2)}(x+1)}{\Gamma(x+1)}+
\biggl(\frac{\Gamma^{(1)}(x+1)}{\Gamma(x+1)}\biggr)^2
\biggr\}=
$$
$$
=P_N(x)\biggl\{\biggl(\log \frac{N}{x}-\frac{1}{2x}+O(x^{-2})\biggr)^2+
O(x^{-1}\,\log x)\biggr\}.
$$
Finally we obtain,
$$
\sum_{n\in B_N}P^{(2)}_N(n)=\sum_{n\in B_N}P_N(n)
\biggl\{\biggl(\log \frac{N}{n}\biggr)^2-
\frac{\log( Nn^{-1})}{n}+O(n^{-1}\,\log n)\biggr\},
$$
where $\widetilde B_N=[Ne^{-\sqrt{\frac{2\log N}{N}}},
Ne^{\sqrt{\frac{2\log N}{N}}}]$.
Note that for $n\in \widetilde B_N $, we have
$|\log^2\frac{N}{n}|\le \frac{\log N}{N}$, where $\log N=O(1)$.
By similar reasons, $\frac{\log (Nn^{-1})}{n}+O(n^{-1}\,\log n)=
O(N^{-1})$ uniformly on this set; in the sense of
asymptotical expansions, $\widetilde B_N$ is equivalent to
the  previously defined set $B_N$.
Hence,
$$
\sum_{n\in\widetilde  B_N}P^{(2)}_N(n)=O(N^{-1})
\sum_{n\in\widetilde  B_N}P_N(n)
=O(N^{-1}).
\eqno(A.5)
$$

\subsection{Corrections to the Gaussian distribution}
From the Stirling expansion (A.2) we have
$$
P_N(n)=\frac{1}{\sqrt{ 2\pi N}}e^{-N+n-n\log \frac{n}{N}}
\frac{1}{\sqrt{ n/N}(1+(12n)^{-1}+O(n^{-2}))}.
$$
Set $N=x\sqrt{N}+N$; then for $x\in D_N$ we obtain
$$
-N+n-n\log \frac{n}{N}=-N+x\sqrt{N}+N-(x\sqrt{N}+N)\log(1+x/\sqrt{N})=
$$
$$
=x\sqrt{N}-(x\sqrt{N}+N)
\biggl\{\frac{x}{\sqrt{N}}-
\frac{x^2}{2N}+\frac{x^3}{3\sqrt{N^3}}+\frac{x^4}{4{N^2}}+
O\biggl(\frac{x^5}{\sqrt{N^5}}\biggr)\biggr\}=
$$
$$
=-\frac{x^2}{2} +\frac{x^3}{6\sqrt{N}} -\frac{x^4}{12{N}}+
O\biggl(\frac{x^5}{\sqrt{N^5}}\biggr).
\eqno(A.6)
$$
This expansion contains corrections to the Gaussian exponent
as in (1.1).
Similarly,
$$
\sqrt{n/N}\biggl(1+(12n)^{-1}+O(n^{-2})\biggr)
=\sqrt{1+\frac{x}{\sqrt{N}}}\biggl(1+\frac{1}{12N}
\biggl( 1-\frac{x}{\sqrt{N}} \biggr)+O(N^{-1})\biggr)=
$$
$$
=\biggl(1+\frac{x}{2\sqrt{N}}\biggr)\biggl(1+\frac{1}{12N}
\biggl( 1-\frac{x}{\sqrt{N}} \biggr)+O(x/N^{3/2})+O(N^{-1})\biggr)=
$$
$$
=1+\frac{x}{2\sqrt{N}}
+O(N^{-1})+O(x/N^{3/2}).
\eqno(A.7)
$$
Equations (A.6)--(A.7) imply expansion (1.1).

\subsection{Komatsu inequalities}
The inequalities of the Komatsu type can easily be justified by the following
construction. For a positive absolutely integrable function $f$,
we consider the integral
$\int_x^{\infty}\,f(y)\, dy=F(x) $, $F(\infty)=0$. If there exists some smooth
positive decreasing function $g(x)$ such that $g(\infty)=0$ and
$$
\frac{d F(x)}{dx}= -f(x)\ge \frac{d g(x)}{dx},
\eqno(A.8)
$$
then $F(x)\le g(x)$. Indeed, the inequality (A.8) for derivatives implies
that $g$ decreases at infinity more rapidly then $F$. Since they
coincide at infinity, this proves that $F(x)\le g(x)$.
In the Komatsu inequality and  in our generalization (3.3),
the right-hand sides satisfy (A.8). For the left-hand sides,
the proof is similar;
it uses the characteristic inequality opposite to (A.8).
}

\end {document}